\documentclass[letterpaper,12pt,peerreviewca,draftcls]{IEEEtran} 

\usepackage{Resources/csm16} 
\usepackage[margin=1in]{geometry}
\usepackage{amsmath} 
\usepackage[compress]{cite}

\usepackage{amsfonts}

\usepackage{amsthm}

\usepackage{comment}

\usepackage{booktabs}

\usepackage{graphicx} 

\usepackage{url}
\usepackage{graphicx,xcolor}
\usepackage{verbatim}
\makeatletter
\newcommand{\verbatimfont}[1]{\def\verbatim@font{#1}}%
\makeatother
\verbatimfont{\ttfamily\small}

\newcommand{\bi}{\begin{itemize}}\newcommand{\ei}{\end{itemize}}
\newcommand{\be}{\begin{equation}}\newcommand{\ee}{\end{equation}}
\newcommand{\bee}{\begin{enumerate}}\newcommand{\eee}{\end{enumerate}}
\newcommand{\bea}{\begin{eqnarray}}\newcommand{\eea}{\end{eqnarray}}
\newcommand{\beas}{\begin{eqnarray*}}\newcommand{\eeas}{\end{eqnarray*}}
\newcommand{\bc}{\begin{center}}\newcommand{\ec}{\end{center}}

\usepackage[left,pagewise]{lineno}
\usepackage[english]{babel} 
\usepackage{blindtext}

\usepackage{todonotes}

\usepackage [autostyle, english = american]{csquotes}
\MakeOuterQuote{"}

\theoremstyle{definition}

\theoremstyle{plain}

\usepackage{bigints}

\title{A Research and Educational Robotic Testbed for Real-time Control of Emerging Mobility Systems: From Theory to Scaled Experiments}
\author{Behdad Chalaki, Logan E. Beaver, A M Ishtiaque Mahbub, Heeseung Bang,\\
 and Andreas A. Malikopoulos\\
	POC: A. A. Malikopoulos (andreas@udel.edu)\\ \today }

\makeatletter
\let\NAT@parse\undefined
\makeatother

\newif\ifPDF \ifx\pdfoutput\undefined\PDFfalse \else\ifnum\pdfoutput > 0\PDFtrue \else\PDFfalse \fi \fi
\ifPDF 
\usepackage[pdftex, plainpages = false, colorlinks=true, linkcolor=green!30!blue, citecolor = green!50!blue, urlcolor = blue, filecolor=black, pagebackref=false, hypertexnames=false,  pdfpagelabels ]{hyperref}
\fi

\usepackage[framemethod=TikZ]{mdframed}
\usepackage{adjustbox}
\usepackage[most]{tcolorbox}
\tcbset{enhanced jigsaw,colback=blue!5!white,colframe=blue!55!green,
fonttitle=\bfseries}

\begin{document}
\maketitle
\CSMsetup 

\begin{tcolorbox}[parbox=false,enhanced,breakable]

\section[Summary]{Sidebar: Summary (abstract)}\label{sb:summary}
Emerging mobility systems, e.g., connected and automated vehicles (CAVs), shared mobility, and electric vehicles, provide the most intriguing opportunity for enabling users to better monitor transportation network conditions and make better decisions for improving safety and transportation efficiency. 
However, before connectivity and automation are deployed en masse, a thorough evaluation of CAVs is required\textemdash ranging from numerical simulation to real-world public roads. 
Assessment of the performance of CAVs in scaled testbeds has recently gained momentum due to the flexibility they offer to conduct quick, repeatable experiments that could go one step beyond simulation.
This article introduces the Information and Decision Science Lab's Scaled Smart City (IDS$^3$C), a $1$:$25$ research and educational scaled robotic testbed that is capable of replicating different real-world urban traffic scenarios. 
IDS$^3$C was designed to investigate the effect of
emerging mobility systems on safety and transportation efficiency.
On the educational front, IDS$^3$C can be used for (a) training and educating graduate students by exposing them to a balanced mix of theory and practice, (b) integrating the research outcomes into existing courses, (c) involving undergraduate students in research, (d) creating interactive educational demos, and (e) reaching out to high-school students.  IDS$^3$C has become a research and educational catalyst for motivating interest in undergraduate and high-school students in science, technology, engineering, and mathematics.
In our exposition, we also present a real-time control framework that can be used to coordinate CAVs at traffic scenarios such as crossing signal-free intersections, merging at roadways and roundabouts, cruising in congested traffic, passing through speed reduction zones, and lane-merging or passing maneuvers. Finally, we provide a tutorial for applying our framework in coordinating robotic CAVs to a multi-lane roundabout scenario and a transportation corridor in IDS$^3$C.

\end{tcolorbox}

Emerging mobility systems, e.g., connected and automated vehicles (CAVs), shared mobility, and electric vehicles, mark a paradigm shift in which a myriad of opportunities exist for users to better monitor the transportation network conditions and make optimal operating decisions to improve safety and reduce pollution, energy consumption, and travel delays \cite{zhao2019enhanced}. 
As we move to increasingly complex emerging mobility systems, new control approaches are needed to optimize the impact on system behavior \cite{Malikopoulos2016c} of the interplay between vehicles at different traffic scenarios \cite{Spieser2014}.  
Several studies have shown the benefits of CAVs to reduce energy consumption and alleviate traffic congestion in specific transportation scenarios \cite{Malikopoulos2018d, Zhao2018CTA, cassandras2019b}.
There have been two major approaches to utilizing connectivity and automation of vehicles, namely, platooning and traffic smoothing. A platoon is defined as a group of closely-coupled vehicles traveling to reduce their aerodynamic drag, especially at high cruising speeds. The concept of platoon formation is a popular system-level approach to address traffic congestion, which gained momentum in the $1980$s and $1990$s \cite{Shladover1991,varaiya1993smart,Rajamani2000}. There has been a rich body of research exploring various methods of forming and/or utilizing platoons to improve transportation efficiency \cite{van2017fuel,wang2017developing,johansson2018multi,karbalaieali2019dynamic,yao2019managing,xiong2019analysis,pourmohammad2020platform,ard2020optimizing,Kumaravel:2021uk,Kumaravel:2021wi,mahbub2021_platoonMixed,bhoopalam2018planning,mahbub2022ACC}. 
Traffic smoothing is another approach that has been explored to mitigate the speed variation of individual vehicles throughout the transportation network, which may be introduced by unnecessary braking and the topology of the road network. 
One of the very early efforts in this direction was proposed by Athans \cite{Athans1969} for safe and efficient coordination of merging maneuvers with the intention of avoiding congestion. Assuming a given merging sequence, Levine and Athans formulated the merging problem as a linear optimal regulator \cite{Levine1966} to control a single string of vehicles, with the aim of minimizing the speed errors that will affect the desired headway between each consecutive pair of vehicles.  
Since then, several studies have been reported in the literature that investigate traffic smoothing to eliminate stop-and-go driving at traffic scenarios, such as single intersections \cite{Dresner2008,makarem2012fluent,Lee2012,Kim2014,azimi2014stip,wu2014distributed,Colombo2014,kamal2014vehicle,Gregoire2014a,Malikopoulos2017,fayazi2018mixed,hult2018optimal,wei2018intersection,bichiou2018developing,mirheli2019consensus,kloock2019distributed,zhang2019decentralized,Malikopoulos2019CDC,malikopoulos2019ACC,tian2020game,pan2020optimal,xu2021comparison,Malikopoulos2020,chalaki2020hysteretic,xu2019cooperative,guney2020scheduling,Dresner2004, DeLaFortelle2010,hausknecht2011autonomous, Huang2012,jin2012multi,levin2016optimizing}, multiple adjacent intersections \cite{hausknecht2011autonomous,du2018hierarchical,ashtiani2018multi,yu2019corridor,Mahbub2019ACC,chalaki2020TITS,chalaki2020TCST,chalaki2019optimal,rodriguez2021distributed,pei2021distributed, Zhao2019CCTA-1}, merging roadways \cite{Rios-Torres2017,Ntousakis2016aa,belletti2017expert,ito2018coordination,Malikopoulos2018b,jing2019cooperative,ding2019rule,xiao2020decentralized,liao2021cooperative}, roundabouts \cite{Zhao2018CTA  ,zohdy2013enhancing,debada2016autonomous,bakibillah2019optimal,chalaki2020experimental,xu2021decentralized,jang2019simulation,chalaki2020ICCA}, speed reduction zones and lane drops \cite{ramezani2015optimized,ma2016freeway,malikopoulos2018optimal,piacentini2019highway,abdelghaffar2020development,goulet2021distributed,nguyen2021system,kreidieh2018dissipating,vinitsky2018lagrangian}, and transportation corridors\cite{xia2013dynamic,roncoli2015traffic,roncoli2016hierarchical,Zhao2018ITSC,mahbub2020ACC-2,mahbub2020decentralized, mahbub2020sae-2}.
Two recent survey papers \cite{Malikopoulos2016a,guanetti2018control} provide a comprehensive review of the state-of-the-art methods and challenges in this area.

Commercial simulation platforms are currently available for testing and validating control algorithms for CAVs in a safe and cost-efficient setting. Simulation can help us gather key information about how the system performs in an idealized environment.
However, evaluating the performance of CAVs in a simulation environment imposes  limitations since modeling the exact vehicle dynamics and driving behavior is not feasible. 
Capturing the complexities arising from data loss and transmission latency associated with connectivity and communication networks can be also challenging. 
As Grim et al. \cite{grim2013simulations} stated, "\textit{the problem with simulations is that they are doomed to succeed.}" Although there have been several studies reporting on the impact of coordination of CAVs in traffic scenarios, e.g., intersections, merging at roadways and roundabouts, the effectiveness of these approaches have been mostly shown in simulation. Therefore, validating control approaches for CAVs in a physical testbed is of great importance.

Scaled testbeds for 
CAVs have attracted considerable attention over the last few years. Such testbeds can be used to conduct quick and repeatable experiments in an effort to go one step beyond simulation. Gulliver \cite{pahlavan2011gulliver} and MOPED \cite{axelsson2014moped} have been the outcome of early efforts on developing scaled testbeds for robotic vehicles. Gulliver's focus is mainly on communication among vehicles, while MOPED is focused on low-level control of a single scaled vehicle. MIT's Duckietown \cite{Paull2017} employs differential drive robots and Go-CHART \cite{9341770} uses four-wheel skid-steer vehicles. Both testbeds focus primarily on local perception and autonomy. The Cambridge Minicars \cite{Hyldmar2019} is another testbed for emulating cooperative driving in highway traffic conditions. A general-purpose robotic testbed called the Robotarium has been developed \cite{wilson2020robotarium}, which features differential drive robots. The Cyber-Physical Mobility Lab \cite{kloock2020cyber} has implemented another scaled testbed on decision-making policies and trajectory planning. For a relatively recent review of such robotics testbeds, see \cite{jimenez2013testbeds}. 

In 2017, we established the Information and Decision Science Lab's Scaled Smart City (IDS$^3$C) to develop and validate control algorithms for emerging mobility systems. IDS$^3$C  occupies a 20 by 20 feet area. It includes 50 robotic cars and 10 aerial vehicles, and can replicate real-world traffic scenarios in a small and controlled environment. This testbed can help us prove concepts beyond the simulation level and understand the implications of errors/delays in the vehicle-to-vehicle and vehicle-to-infrastructure communication as well as their impact on energy usage. IDS$^3$C can help us implement control algorithms for coordinating CAVs in different traffic scenarios, such as intersections, merging roadways, speed reduction zones, roundabouts, and transportation corridors.  IDS$^3$C also includes driver emulation stations interfaced directly with the cars which allow us to explore human driving behavior. The robotic cars share many features with full-scale cars, such as four wheels, built-in suspension, and an Ackermann steering mechanism.

There are several features that distinguish IDS$^3$C from other testbeds. 
First, unlike MIT's Duckietown \cite{Paull2017} and GO-Chart \cite{9341770}, the CAVs in IDS$^3$C resemble full-scale vehicles by using four wheels, built-in suspension, and an Ackermann steering mechanism. 
Second, in contrast to the scaled testbeds reported in  \cite{Paull2017,Hyldmar2019,wilson2020robotarium,kloock2020cyber}, IDS$^3$C is equipped with driver emulation stations that interface directly with the robotic cars. These stations enable us to explore and study human driving behavior and their interactions with CAVs. 
Being able to study how CAVs can safely interact and co-exist with human-driven vehicles is of great importance since different penetration rates of CAVs can significantly alter transportation efficiency and safety. 
Third, rather than focusing on specific scenarios in a transportation network \cite{pahlavan2011gulliver,Hyldmar2019} or a single individual vehicle \cite{axelsson2014moped}, IDS$^3$C can accommodate almost every possible traffic scenario, including crossing three- and four-way intersections, merging at roadways and roundabouts, cruising in congested traffic, passing through speed reduction zones, and lane-merging or passing maneuvers. These features make  IDS$^3$C a unique scaled robotic testbed to study problems in emerging mobility systems such as coordination of CAVs, shared-mobility, eco-routing, and firs/last-mile delivery.  Finally, only a few testbeds \cite{kloock2020cyber,Paull2017} are equipped with a ``digital twin.'' The digital twin of IDS$^3$C, called the Information and Decision Science Lab's Scaled Smart Digital City (IDS $3$D City), is a Unity-based virtual simulation environment that can operate alongside the physical IDS$^3$C and interface with the existing control framework. The IDS $3$D City provides the framework to develop and implement control algorithms for emerging mobility systems in simulation before moving to the physical IDS$^3$C for validation. More details about the IDS $3$D City can be found in \cite{Ray2021DigitalCity}.

In what follows, we start our exposition by providing a brief description of the hardware and software architecture of IDS$^3$C.
Then, we present an overview of a real-time coordination framework for CAVs which has been implemented and validated in IDS$^3$C and field testing \cite{mahbub2020sae-1,mahbub2020sae-2}. Next, we present a tutorial of this framework in an application to a multi-lane roundabout in IDS$^3$C using $9$ CAVs. 
Finally, we  provide a demonstration study in IDS$^3$C by using a fleet of $15$ CAVs and show how we can improve traffic throughput along a transportation corridor, which consists of a roundabout, an intersection, and a merging roadway.

\section{Information and Decision Science Lab's Scaled Smart City}

\begin{figure}[htbp]
    \centering
    \includegraphics[width=0.85\textwidth]{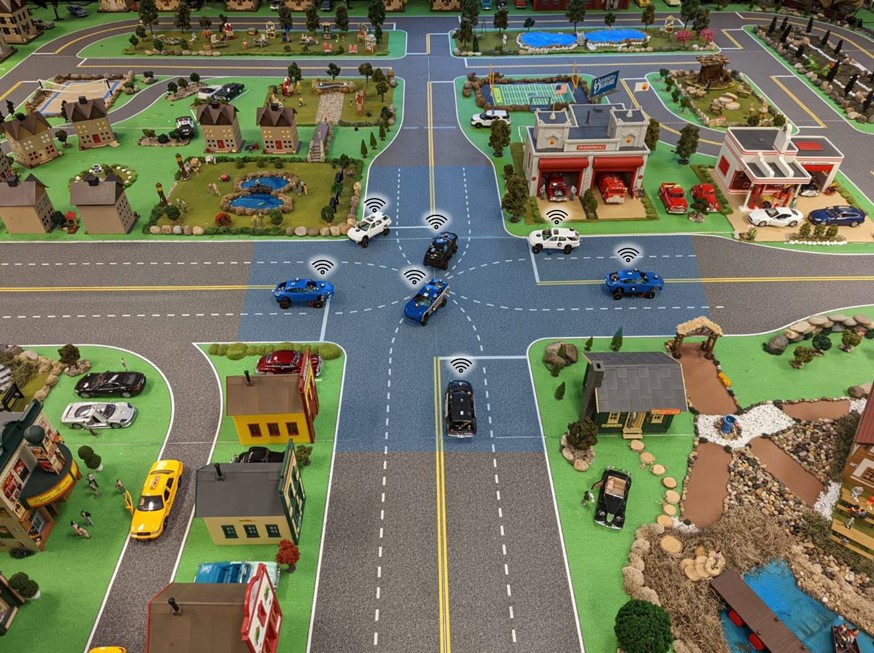}
    \caption{A view of the IDS Lab's Scaled Smart City (IDS$^3$C) with connected and automated vehicles coordinating at an intersection.}
    \label{fig:IDS3C}
\end{figure}

IDS$^3$C (Fig. \ref{fig:IDS3C}) is a $1$:$25$ scaled robotic testbed spanning over $400$ square feet, and it is capable of replicating real-world traffic scenarios in a small and controlled environment using $50$ ground and $10$ aerial vehicles.
IDS$^3$C provides the opportunity to prove concepts beyond simulation and understand the implications of errors and delays in the vehicle-to-vehicle and vehicle-to-infrastructure communication as well as their impact on energy usage. IDS$^3$C can also be used to understand the implications of emerging mobility systems, consisting of CAVs, shared mobility, and electric vehicles, on safety and transportation efficiency. Another facet of research that can be explored using IDS$^3$C is complex missions that include the cooperation of aerial and ground vehicles for logistic problems, such as last-mile delivery. 
IDS$^3$C includes driver emulation stations \cite{driver} interfaced directly with the robotic cars which allow us to explore human driving behavior.

IDS$^3$C is equipped with a VICON motion capture system and uses eight cameras to track the position of each vehicle with sub-millimeter accuracy.
The testbed contains a dozen traffic bottlenecks, including merging roadways, multi-lane roundabouts, adjacent intersections, multi-lane intersections, lane-drops, and speed reduction zones. A central mainframe computer (Processor: Intel Core i$7$-$6950$X CPU @ $3.00$ GHz x $20$, Memory: $128$ GB) stores a map of the IDS$^3$C as a database of line and arc segments that make up the road network.
Coordination of the CAVs within the IDS$^3$C is achieved using a multi-level control framework spanning the mainframe computer and the individual CAVs in an experiment.
Each CAV is given its own thread on the central mainframe computer. The latter communicates the vehicle's position through VICON and generates its trajectory.
Lane and reference trajectory tracking are accomplished onboard each CAV in a purely distributed manner.

We developed IDS$^3$C with the capacity to experimentally validate a wide variety of urban mobility scenarios.
This includes eco-routing, mixed traffic \cite{Zhao2018CTA}, ride sharing \cite{zhao2019}, last-mile delivery \cite{K122021AExperiment}, and air-ground coordination \cite{remer2019multi}.
In several recent efforts, we have used IDS$^3$C to implement and validate control algorithms for coordinating CAVs at traffic scenarios, such as merging roadways \cite{Malikopoulos2018b}, roundabouts \cite{chalaki2020experimental}, intersections  \cite{Malikopoulos2019CDC}, adjacent intersections \cite{chalaki2020TCST,chalaki2020TITS,Mahbub2019ACC}, and corridors \cite{Beaver2020DemonstrationCity}. 
We have also used IDS$^3$C to transfer policies from neural networks \cite{jang2019simulation,chalaki2020ICCA} and handle the stochasticity that arises in physical systems \cite{chalaki2021RobustGP}. IDS$^3$C is in the position to provide a means for user interaction through a mobile application, which enables submitting origin-destination travel requests for dynamic routing in shared mobility and last-mile delivery scenarios. More recently, we introduced a Unity-based virtual simulation environment for emerging mobility systems, called the IDS $3$D City, intended to operate alongside its physical peer, IDS$^3$C, and interface with its existing control framework. For a brief summary of IDS $3$D City see "\nameref{sb:IDS 3D City}", and for further technical discussions discussion see \cite{Ray2021DigitalCity}.

\begin{tcolorbox}[parbox=false,enhanced,breakable]
\section[The Information and Decision Science Lab's Scaled Smart Digital City (IDS 3D City)]{{Sidebar: Information and Decision Science Lab's Scaled Smart Digital City (IDS 3D City)}} \label{sb:IDS 3D City}
The Information and Decision Science Lab's Scaled Smart Digital City (IDS 3D City) is a digital replica of the Information and Decision Science Lab's Scaled Smart City (IDS$^3$C) using AirSim and Unity.
We have designed the IDS 3D City to integrate the control framework used in IDS$^3$C to simulate virtual vehicles.
The IDS 3D City enables users to rapidly iterate their control algorithms and experiment parameters before deploying it to IDS$^3$C.
A schematic of how the IDS 3D City interacts with IDS$^3$C
is shown in Fig. \ref{fig:commGraph}. The end result is a transition between the physical \text{and} virtual environments with minimal changes to input files, as well as the capability to mix physical and virtual vehicles.
\par
    {\centering
    \includegraphics[width=1\textwidth]{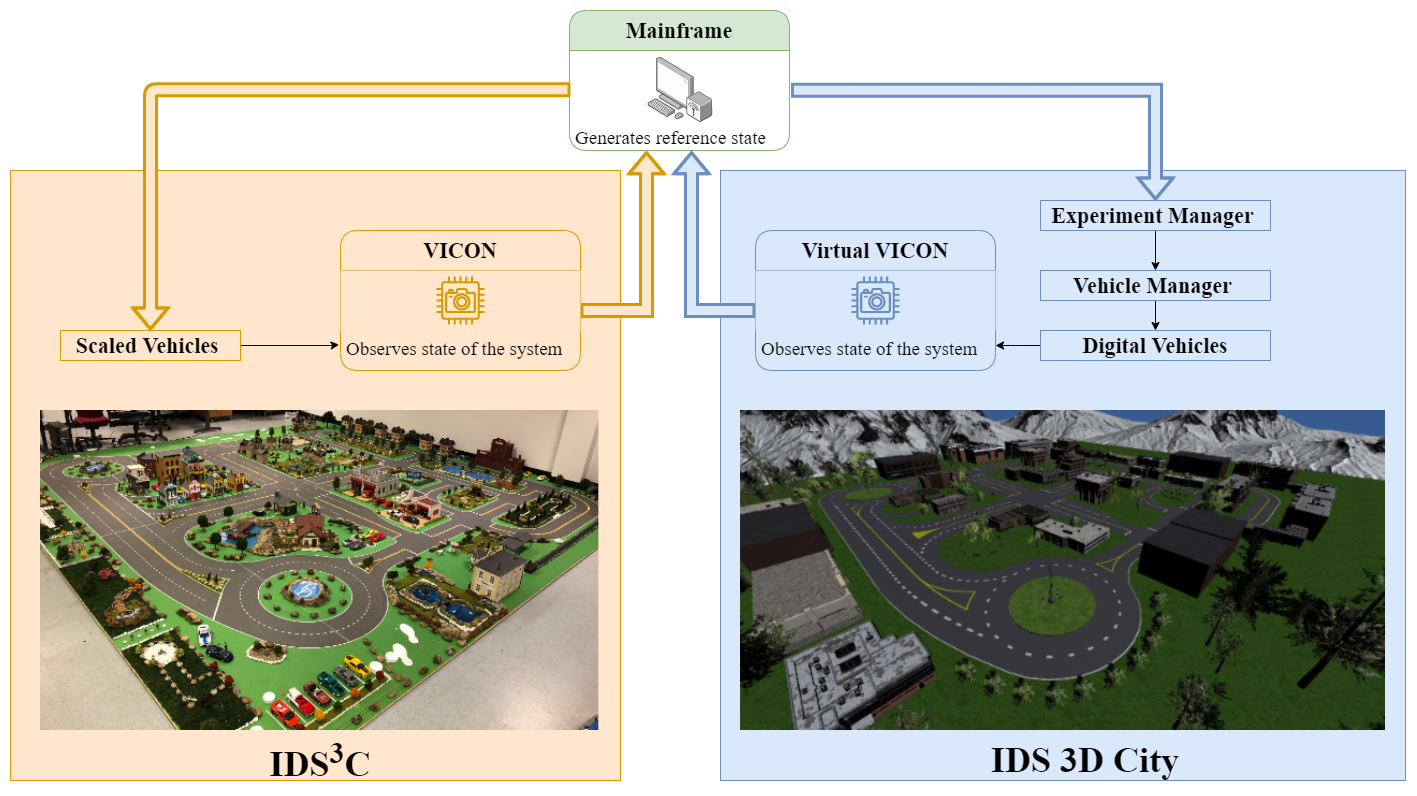}
    \captionof{figure}{The physical and virtual city environments. The mainframe computer can switch between physical and virtual experiment.}
    \label{fig:commGraph}}
\par
\end{tcolorbox}

We have used IDS$^3$C to develop a control framework for  real-time coordination of CAVs at different traffic scenarios such as crossing signal-free intersections, merging at roadways and roundabouts, cruising in congested traffic, passing through speed reduction zones, and lane-merging or passing maneuvers.
Next, we outline the main features of this framework.

\section{A Control Framework for Coordination of Connected and Automated Vehicles}

IDS$^3$C has been extensively used for the development and implementation of control algorithms aimed at coordinating CAVs in different traffic scenarios.
By coordinating CAVs in traffic scenarios, they do not have to come to a full stop, thereby conserving momentum and fuel while also improving travel time. 
In this context, safety is explicitly guaranteed by imposing constraints on each vehicle, including rear-end safety, maximum speed limit, and lateral collision avoidance.

Several research efforts have considered a two-level optimization framework for coordinating CAVs consisting of a travel time minimization (upper-level) and an energy minimization problem (lower-level).
For each CAV with a given origin-destination, the solution of the upper-level problem yields the travel time for the CAV to exit a ``control zone,'' inside of which the CAVs can communicate with each other.
The solution of the low-level problem yields, for each CAV, the control input (acceleration/deceleration) to achieve the solution of the upper-level problem while minimizing energy consumption subject to the state, control, and safety constraints. 
The details of the low-level problem can be found in \cite{Malikopoulos2017,Malikopoulos2020, mahbub2020Automatica-2}. Solving a constrained optimal control problem leads to a system of nonlinear equations that are often infeasible to solve in real time.  
For more details on the constrained optimal control and the associated technical challenges, see "\nameref{sb:constrained optimal}." 

\begin{tcolorbox}[parbox=false,enhanced,breakable]
\section[The Challenge with Constrained Optimal Control]{{Sidebar: The Challenge with Constrained Optimal Control}} \label{sb:constrained optimal}

The standard methodology to solve a continuous-time constrained optimal control problem is to employ Hamiltonian analysis with interior point state and/or control constraints \cite{S1}. Namely, we first start with the unconstrained arc and derive the solution of the optimal control problem without considering any of the state or control constraints. 
If the unconstrained solution violates any of the state or control constraints, then the unconstrained arc is pieced together with the arc corresponding to the violated constraint. The two arcs yield a set of algebraic equations that are solved simultaneously using the boundary conditions and optimality conditions between the arcs. If the resulting solution, which includes the determination of the optimal
switching time from one arc to the next one, violates another constraint, then the last two arcs are pieced together with the arc corresponding to the new violated constraint, and we re-solve the problem with the three arcs pieced together. The three arcs will yield a new set of algebraic equations that need to be solved, and this process is repeated until the solution does not violate any constraints.  
This iterative process can be computationally intensive for several reasons. First, the Euler-Lagrange equations are numerically unstable for non-conservative systems, leading to significant numerical challenges \cite{S2}.
Second, the number of active constraints is not known a priori, and it may require a significant number of iterations to compute.
Third, the boundary conditions and recursive equations may be implicit functions that do not have a closed-form analytical solution.

Excluding cases with terminal speed and safety constraints, in recent work \cite{S3,S4}, we have introduced a condition-based solution framework for the optimal coordination of CAVs, which leads to a closed-form analytical solution without this iterative procedure. In this framework, we mathematically characterized the activation cases of different state and control constraint combinations, and provided a set of a priori conditions under which different constraint combinations can become active. Although this approach alleviates the computational complexity of the constrained optimal control in the coordination problem to some extent, the aforementioned iterative procedure is still required for cases when safety and terminal speed constraints are included.

\end{tcolorbox}

To avoid the challenges associated with the constrained optimal control, we have proposed an alternative control framework consisting of a single-level optimization aimed at both minimizing energy consumption and improving traffic throughput \cite{Malikopoulos2020}. Next, we highlight the features of this framework.

\subsection{Problem Formulation}
\begin{figure}[htbp]
	\centering
	\includegraphics[width=0.50\textwidth]{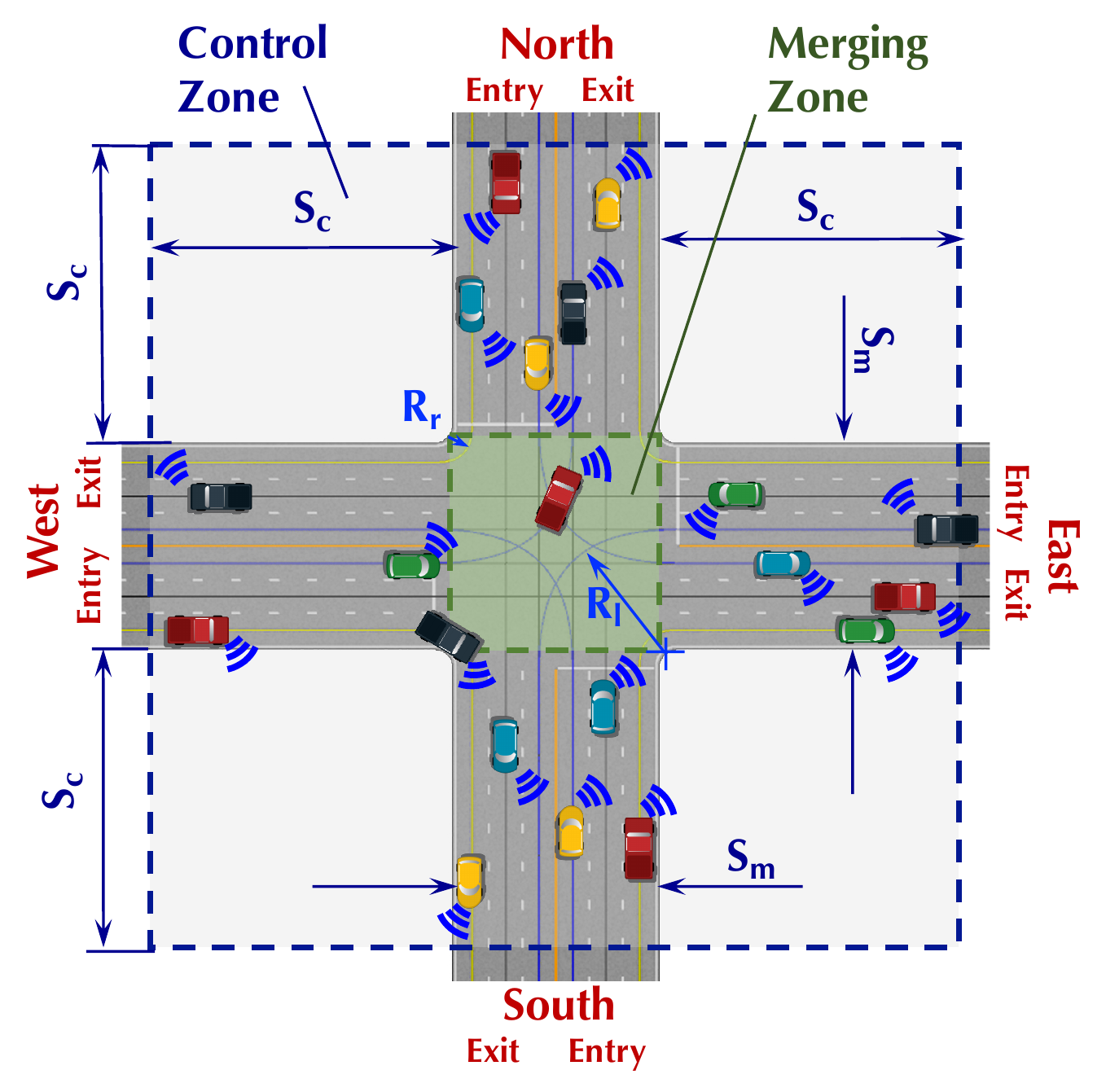} 
	\caption{A signal-free intersection with connected and automated vehicles.}
	\label{fig:1}
\end{figure}

Although the control framework presented here can be applied to any traffic scenario, we use an intersection as a reference to provide the fundamental ideas. This is because an intersection has unique features which makes it technically more challenging compared to other traffic scenarios. However, our analysis can be applied to other traffic scenarios.
We consider CAVs at a $100$\% penetration rate crossing a signal-free intersection (Fig. \ref{fig:1}). The region at the center of the intersection, called the \textit{merging zone}, is the area of potential lateral collision of  CAVs. The intersection has a \textit{control zone} inside of which the CAVs plan their time trajectories (a time trajectory yields the time that a CAV is on a given position inside the control zone) 
by  communicating with each other and with a \textit{coordinator}, i.e., a roadside unit that stores the planned time trajectories of each CAV as they pass through the control zone. The distance from the entry of  the control zone until the entry of the merging zone is $S_c$. Although it is not restrictive, we consider $S_c$ to be the same for all entry points of the control zone. We also consider the merging zone to be a square of side $S_m$ (Fig. \ref{fig:1}). Note that  $S_c$ could be in the order of hundreds of meters depending on the CAVs' communication range capability, while $S_m$ is the length of a typical intersection. The CAVs crossing the intersection can also make a right turn of radius $R_r$, or a left turn of radius $R_l$ (Fig. \ref{fig:1}). The aforementioned values of the intersection's geometry are not restrictive in our modeling framework, and are used only to determine the total distance traveled by each CAV inside the control zone.
In our problem formulation, we assume that each CAV can communicate with other CAVs and the coordinator without any errors or delays. It is relatively straightforward to relax this assumption as long as the noise in the communication, measurements, and/or delays is bounded. We also assume that upon entering the control zone, the initial state of each CAV is feasible, that is, none of the speed or safety constraints are violated. This is a reasonable assumption since CAVs are automated; therefore, there is no compelling reason for them to violate any of the constraints by the time they enter the control zone.

We denote the set of CAVs in the control zone by the set $\mathcal{N}(t) = \{1, \ldots, N(t)\}$, where $N(t)\in\mathbb{N}$ is the total number of CAVs at time $t\in\mathbb{R}_{\geq 0}$.
In our framework, for each CAV $i\in\mathcal{N}(t)$, we seek to jointly minimize energy consumption and travel time.
Upon entering the control zone, CAV $i$ communicates with the coordinator and receives the time trajectory of all CAVs $j\in\mathcal{N}(t)\setminus\{i\}$.
Next, CAV $i$ computes the time $t_i^f$ that it must exit the control zone while guaranteeing that its corresponding energy optimal time trajectory does not activate any of the state, control, and safety constraints. This trajectory is communicated back to the coordinator for the subsequent CAVs to plan their trajectories with the same sequence as they enter the control zone. If two or more CAVs enter the control zone simultaneously, then the coordinator arbitrarily determines the sequence in which they receive information to plan their trajectories.
Addressing the problem sequentially makes coordination between the CAVs tractable at the possible cost of selecting a sub-optimal planning sequence.
Finding the optimal sequence of decision making is a combinatorial problem, which is NP-hard \cite{de2013autonomous}.
In a recent paper, we reported how vehicles can dynamically change their decision-making sequence and replan to improve the throughput at an intersection \cite{chalaki2021Reseq} by relaxing the first-come-first-serve decision making.

By enforcing the unconstrained energy optimal time trajectory that guarantees that none of the state, control, and safety constraints becomes active, we avoid the challenges associated with the real-time implementation of the constrained optimal control solution.
In our analysis, we consider that each CAV $i\in\mathcal{N}(t)$ is governed by the following dynamics,
\begin{align}\label{eq:model2}
\dot{p}_{i}(t)  =v_{i}(t), \quad \dot{v}_{i}(t)  =u_{i}(t), \quad \dot{s}_{i}(t) = v_{k}(t)-v_{i}(t), ~t\in[t_i^0, t_i^f],
\end{align}
where $t_i^0$ and $t_i^f$ correspond to the times that CAV $i$ enters and exits the control zone, respectively; $p_{i}(t)\in\mathcal{P}_{i}$ is the position of each CAV $i$ from the entry until the exit of the control zone; $v_{i}(t)\in\mathcal{V}_{i}$ and $u_{i}(t)\in\mathcal{U}_{i}$ are the speed and acceleration/deceleration (control input), respectively, of each CAV $i$ inside the control zone; ~$s_{i}(t)\in\mathcal{S}_{i}$ denotes the distance of CAV $i$ from CAV $k$ which is physically located  ahead of $i$, and $v_k(t)$ is the speed of CAV $k$. The sets $\mathcal{P}_{i}, \mathcal{V}_{i}, ~\mathcal{U}_{i}$, and $\mathcal{S}_{i}$, $i\in\mathcal{N}(t)$, are complete and totally bounded subsets of $\mathbb{R}$.

In our framework, we impose the following constraints to ensure the CAV's control input and state remain within an admissible range,
\begin{align}\label{eq:uvConstraint}
    u_{i,\min} \leq u_i(t) \leq u_{i,\max}, \quad
    0 < v_{\min} \leq v_i(t) \leq v_{\max},
\end{align}
for all $t\in[t_i^0, t_i^f]$, where $u_{i,\min},u_{i,\max}$ are the minimum and maximum control inputs and $v_{\min},v_{\max}$ are the minimum and maximum speed limit, respectively.
To ensure that no rear-end collisions occur between two CAVs traveling in the same lane, we impose the rear-end safety constraint,
\begin{equation} \label{eq:rearSafety}
    s_i(t) = p_k(t) - \lambda_k - p_i(t) \geq \delta_i(t) = \gamma + \rho_i v_i(t),
\end{equation}
where $\lambda_k$ is the length of CAV $k$, $\gamma$ is the standstill distance, and $\rho_i$ is the minimum time headway that CAV $i$ wishes to maintain with preceding CAV $k$.

Finally, let $j\in\mathcal{N}(t) \setminus \{i\}$ correspond to another CAV that has already entered the control zone and may have a lateral collision with CAV $i$. For example, suppose CAV $i$ travels north-south and CAV $j$ travels east-west (see Fig. \ref{fig:1}).
Then there is a \textit{conflict point} where the paths of $i$ and $j$ intersect, and hence, a potential lateral collision might occur. We include all such conflict points in a finite set $\mathcal{O}\subset\mathbb{N}$, which is entirely determined by the geometry of the roads.
Let $p_{i}^n$ and $p_{j}^n$ be the position of the conflict point $n\in\mathcal{O}$ along the paths of CAV $i$ and $j$, respectively.
CAV $i$ can cross this conflict point either after, or before CAV $j$. 
In the first case, we have
\begin{equation} \label{eq:lateralBefore}
    p_i^n - p_i(t) \geq \delta_i(t), \quad \text{for all}\quad t\in[t_i^0, t_j^n],
\end{equation}
where $t_j^n$ is the known time that CAV $j$ reaches at conflict point $n$, that is, position $p_j^n$.
In the second case, we have
\begin{equation} \label{eq:lateralAfter}
    p_j^n - p_j(t) \geq \delta_j(t), \quad \text{for all}\quad t\in[t_j^0, t_i^n],
\end{equation}
where $t_i^n$ is determined by the trajectory planned by CAV $i$.

Since $ 0 < v_{\min} \leq v_i(t)$, the position $p_i(t)$ is a strictly increasing function.
Thus, the inverse $t_i\left(\cdot\right) = p_i^{-1}\left(\cdot\right)$ exists and it is called the \textit{time trajectory} of CAV $i$, hence we have $t_i^n = p_i^{-1}\left(p_i^n\right)$. 
The closed-form solution of the inverse function is derived in \cite{Malikopoulos2020}. 
To guarantee lateral safety between CAV $i$ and CAV $j$ at a conflict point $n$, either \eqref{eq:lateralBefore} or \eqref{eq:lateralAfter} must be satisfied. Therefore, we impose the following lateral safety constraint on CAV $i$,
\begin{align}
    \min \Bigg\{ &\max_{t\in[t_i^0, t_j^n]} \{ \delta_i(t) + p_i(t) - p_i^n\}, \max_{t\in[t_j^0, t_i^n]} \{ \delta_j(t) + p_j(t) - p_j^n \}   \Bigg\} \leq 0 . \label{eq:lateralMinSafety}
\end{align}
When CAV $i\in\mathcal{N}(t)$ enters the control zone it must determine the exit time $t_i^f$ such that the resulting time trajectory does not activate any of \eqref{eq:model2} - \eqref{eq:rearSafety} and \eqref{eq:lateralMinSafety}.
The unconstrained solution for CAV $i$ is  
\begin{align} \label{eq:optimalTrajectory}
    u_i(t) = 6 a_i t + 2 b_i, \quad
    v_i(t) = 3 a_i t^2 + 2 b_i t + c_i, \quad
    p_i(t) = a_i t^3 + b_i t^2 + c_i t + d_i, 
\end{align}
where $a_i, b_i, c_i, d_i$ are constants of integration.
CAV $i$ must also satisfy the boundary conditions $ \Big( p_i(t_i^0),\, v_i(t_i^0)   \Big)  = \Big( 0,\, v_i^0 \Big)$ and
    $\Big( p_i(t_i^f),\, u_i(t_i^f)   \Big)  = \Big( S_i,\, 0 \Big)$,
where $S_i$ is the known length of CAV $i$'s path in the control zone.
For the details on deriving the unconstrained solution, see "\nameref{sb:Unconstrained Hamiltonian}." 

\begin{tcolorbox}[parbox=false,enhanced,breakable]
\section[Unconstrained Optimal Control]{Sidebar: Unconstrained Optimal Control and Boundary Conditions}\label{sb:Unconstrained Hamiltonian}

Let $t_i^{f}$ be the specified exit time of CAV $i$ from the control zone. 
To minimize the energy consumption of $i$ inside the control zone, we minimize transient engine operation through the $L^2$-norm of the control input $u_i(t)$ over the interval $[t_i^0,~t_i^f]$, which is known to have direct benefit in fuel consumption and emission in conventional vehicles \cite{S5,S6}. Namely, CAV $i$ minimizes the following cost function,
\begin{equation}\label{eq:cost_func_1}
    J_i(u_i(t),t_i^{f})= {\dfrac{1}{2}} \bigints_{t_i^{0}}^{t_i^{f}} u_{i}(t)^2~dt.
\end{equation}
For each CAV $i$ in the control zone, the unconstrained Hamiltonian is
\begin{equation}\label{c1}
\begin{aligned}
H_i(t,p_i(t),v_i(t),u_i(t))&=\frac{1}{2}u_i(t)^2+\lambda_i^p v_i (t) +\lambda_i^v u_i(t),
\end{aligned}
\end{equation}
where \(\lambda_i^p\) and \(\lambda_i^v\) are costates corresponding to the position and speed of the CAV, respectively.
The Euler-Lagrange optimality equations are
\begin{align}
\dot{\lambda}_i^p=-\frac{\partial H_i}{\partial p_i}&=0, \quad
\dot{\lambda}_i^v=-\frac{\partial H_i}{\partial v_i}= -\lambda_i^p,\label{euler:optimallambdap}\\
\frac{\partial H_i}{\partial u_i}&=u_i+\lambda_i^v=0\label{euler:optimalU}.
\end{align}
Since the speed of CAV $i$ is not specified at the terminal time $t_i^f$, we have \cite{S7}
\begin{equation} \label{eq:lambdaV-sb}
    \lambda_i^v(t_i^f) = 0.
\end{equation}
Applying the Euler-Lagrange optimality conditions \eqref{euler:optimallambdap}-\eqref{euler:optimalU} to the Hamiltonian \eqref{c1} yields $u^\ast_i(t) = -{\lambda_i^v}^\ast= a_i^{\prime}t+b_i^{\prime}$, where $a_i^{\prime}$ and $b_i^{\prime}$ are constants of integration. By integrating the control input, we can find the optimal position and speed trajectories as $p_i(t) =\frac{1}{6}a_i^{\prime}t^3+\frac{1}{2}b_i^{\prime}t^2+c_i^{\prime}t+d_i^{\prime},$ and $v_i(t) =\frac{1}{2}a_i^{\prime}t^2+b_i^{\prime}t+c_i^{\prime},$
where \(a_i^{\prime},b_i^{\prime},c_i^{\prime},d_i^{\prime}\) are constants of integration, which are found by substituting the boundary conditions.
The boundary conditions for any CAV $i$ are $p_i(t_i^0) = p_i^0, ~ v_i(t_i^0) = v_i^0,~ p_i(t_i^f) = p_i^f, \text{ and } u_i(t_i^f) = 0,$
where $p_i$ is known at $t_i^0$ and $t_i^f$ by the geometry of the road, and $v_i^0$ is the speed at which the CAV enters the control zone.
The final boundary condition, $u_i(t_i^f) = 0$, arises from substituting \eqref{eq:lambdaV-sb} into \eqref{euler:optimalU} at $t_i^f$, that is, $ u_i(t_i^f) + \lambda_i^v(t_i^f) = 0$,
which implies $u_i(t_i^f) = 0$.

\end{tcolorbox}

There are five unknown variables that determine the optimal time trajectory of CAV $i$, four constants of integration from \eqref{eq:optimalTrajectory}, and the unknown exit time $t_i^f$.
Without loss of generality, letting $t_i^0 = 0$ implies that
    $p_i(t_i^0) = d_i = 0$, 
    $v_i(t_i^0) = c_i = v_i^0$, 
while $u_i(t_i^f) = 0$ yields
$a_i = \frac{-b_i}{3 t_i^f}$, 
and $p_i(t_i^f) = S_i$ gives
    $b_i = \frac{3(S_i - v_i^0 t_i^f)}{2 (t_i^f)^2}$. 
%
Furthermore,  $t_i^f$ takes a value from a compact set, $[\underline{t}_i^f, \overline{t}_i^f]$. 
See "\nameref{sb:FinalTimeBounds}" for more details on the derivation $\underline{t}_i^f$ and $\overline{t}_i^f$ based on the speed and control constraints and boundary conditions.

\begin{tcolorbox}[parbox=false,enhanced,breakable]
\section[Bounds for Feasible Exit Time]{Sidebar: Derivation of Bounds for Feasible Exit Time}\label{sb:FinalTimeBounds}

The unconstrained optimal trajectory of CAV $i\in\mathcal{N}(t)$ takes the form
\begin{align} 
    u_i(t) &= 6 a_i t + 2 b_i,\label{eqS:UOpt}  \\
    v_i(t) &= 3 a_i t^2 + 2 b_i t + c_i,\label{eqS:VOpt} \\
    p_i(t) &= a_i t^3 + b_i t^2 + c_i t + d_i,\label{eqS:POpt}
\end{align}
where $a_i, b_i, c_i, d_i$ are constants of integration, which are found by using the boundary conditions.
We derive the upper and lower bounds on the exit time of the control zone for a CAV $i\in\mathcal{N}(t)$ using the speed and control constraints by exploiting two properties of the optimal trajectory.
As the optimal control input is linear and satisfies $u_i(t_i^f) = 0$, it must be zero, strictly decreasing, or strictly increasing.
In all three cases $u_i(t)$ achieves its extreme at $t_i^0$, and therefore satisfying $u_{\min} \leq u_i(t_i^0) \leq u_{\max}$ is necessary and sufficient condition to guarantee constraint satisfaction.
Likewise, the speed of CAV $i$ starts at $v_i(t_i^0) = v_i^0\in[v_{\min}, v_{\max}]$ and must be constant, strictly increasing, or strictly decreasing inside the control zone.
In all three cases $v_i(t)$ takes its extreme value at $t_i^f$, and thus satisfying $v_{\min} \leq v_i(t_i^f) \leq v_{\max}$ is necessary and sufficient condition to guarantee constraint satisfaction.

Next, without loss of generality, let $t_i^0=0$ and $p_i^0=0$.
This implies that
    $p_i(t_i^0) = d_i = 0$ and  
    $v_i(t_i^0) = c_i = v_i^0$, 
while $u_i(t_i^f) = 0$ implies
\begin{equation}
a_i = \frac{-b_i}{3 t_i^f}, \label{eqS:ai}
\end{equation}
and $p_i(t_i^f) = S_i$ yields
\begin{equation}
    b_i = \frac{3(S_i - v_i^0 t_i^f)}{2 (t_i^f)^2}. \label{eqS:bi}
\end{equation}
In order to compute the lower bound on exit time of the control zone for CAV $i$, $\underline{t}_{i}^f$, there are two cases to consider:

\textbf{Case L$1$:} CAV $i$ achieves its maximum control input at entry of the control zone, that is, $u_i(t_i^0)=u_{\max}$.
In this case, evaluating \eqref{eqS:UOpt} at $t_i^0=0$ yields
\begin{align}\label{eq:tfumax}
    u_i(t) &= 2 b_i = u_{\max}.
\end{align}
Substituting \eqref{eqS:bi} into \eqref{eq:tfumax} and solving for $t_i^f$ yields the quadratic equation $u_{\max}{t_i^f}^2 + 3v_i^0 t_i^f - 3S_i = 0$,
which has two real roots with opposite signs, since $t_{i,1}^f t_{i,2}^f = \frac{-3 S_i}{u_{\max}}<0$.
Thus, $t_{i,u_{\max}}^f>0$ is $t_{i,u_{\max}}^f = \frac{\sqrt{9 {v_i^0}^2 + 12 S_i u_{\max}} - 3 {v_i^0}}{2 u_{\max}}$.

\textbf{Case L$2$:} CAV $i$ achieves its maximum speed at the end of control zone, that is, $v_i(t_i^f)=v_{\max}$.
For this case, by \eqref{eqS:VOpt}, we have
\begin{equation} \label{eq:vmaxPoly}
    v_i(t_i^f) = {3}a_i {t_i^f}^2 + 2b_i t_i^f + v_i^0 = v_{\max}.
\end{equation}
Substituting \eqref{eqS:ai} and \eqref{eqS:bi} into \eqref{eq:vmaxPoly} yields
\begin{align}
     v_i(t_i^f) &= 3 \Big(\frac{-b_i}{3 t_i^f}\Big) {t_i^f}^2 + 2 b_i t_i^f + {v_i^0}
     = b_i t_i^f + {v_i^0} = \frac{3(S_i - v_i^0 t_i^f)}{2 t_i^f} + {v_i^0} = v_{\max},
\end{align}
which simplifies to $t_{i,v_{\max}}^f = \frac{3 S_i}{v_i^0 + 2 v_{\max}}$.
Thus, our lower bound on $t_i^f$ is given by $\underline{t}_{i}^f = \min \left\{ t_{i,u_{\max}}^f, t_{i,v_{\max}}^f \right \}$.

The upper bound for $t_{i}^f$ can be derived following similar steps for the lower bound, and can be broken into two cases.

\textbf{Case U$1$:} CAV $i$ achieves its minimum control input at the entry of the control zone, that is, $u_i(t_i^0) = u_{\min}$.
This implies $u_{\min}{t_i^f}^2 + 3v_i^0 t_i^f - 3 S_i = 0$,
which has two positive roots, as $t_{i,1}^f t_{i,2}^f = \frac{-3S_i}{u_{\min}}>0$, from which we select the smaller one,
\begin{align} \label{eq:tfUmin}
    t_{i,u_{\min}}^f = \frac{\sqrt{9 {v_i^0}^2 + 12 S_i u_{\min}} -3 v_i^0 }{2 u_{\min} },
\end{align}
as the speed of the vehicle should be always greater than zero.
Note that when $9 {v_i^0}^2 + 12 S_i u_{\min} < 0$ there is no real value of $t_i^f$ which satisfies all of the boundary conditions simultaneously, and therefore the constraint $u(t_i^0) = u_{\min}$ can never become active if \eqref{eq:tfUmin} is complex.
In that case, the upper bound must be given by Case U2.

\textbf{Case U$2$}: CAV $i$ achieves its minimum speed at the entry of the control zone, that is, $v_i(t_i^f) = v_{\min}$.
Evaluating \eqref{eqS:VOpt} at $t_i^f$ yields $v_i(t_i^f) = {3}a_i {t_i^f}^2 + 2b_i t_i^f + v_i^0 = v_{\min}$, in which
substituting \eqref{eqS:ai} and \eqref{eqS:bi} yields
\begin{align}
     v_i(t_i^f) &= 3 \Big(\frac{-b_i}{3 t_i^f}\Big) {t_i^f}^2 + 2 b_i t_i^f + {v_i^0} = b_i t_i^f + {v_i^0} = \frac{3(S_i - v_i^0 t_i^f)}{2 t_i^f} + {v_i^0} = v_{\min},
\end{align} 
which simplifies to $t_{i,v_{\min}}^f = \frac{3 S_i}{v_i^0 + 2 v_{\min}}$.

 Thus, the upper bound on the exit time for CAV $i$ is 
\begin{equation}
        t_{i}^f = 
        \begin{cases}
            t^f_{i,v_{\min}},&\text{ if }\ 9 {v_i^0}^2 + 12 S_i u_{i,\min} < 0,\\
            \max \{ t_{i,u_{\min}}^f  , t_{i,v_{\min}}^f\} ,&\text{otherwise.}\ \\
        \end{cases}
\end{equation}
 where $t^f_{i,v_{\min}} = \frac{3 S_i}{v_i^0 + 2 v_{\min} }$ and $t^f_{i,u_{\min}} = \frac{\sqrt{9 {v_i^0}^2 + 12 S_i u_{\min}} -3 v_i^0 }{2 u_{\min} }$.
\end{tcolorbox}

This leads to the following optimization problem
\begin{align}\label{eq:tif}
    &\min_{t_i^f\in[\underline{t}_i^f, \overline{t}_i^f]} t_i^f \\
    \text{subject to: } & \text{rear-end safety } \eqref{eq:rearSafety}, \text{ lateral safety } \eqref{eq:lateralMinSafety}, \text{ dynamics } \eqref{eq:optimalTrajectory}, \notag \\
    &d_i = 0, \quad c_i = v_i^0, \quad a_i = -\frac{-b_i}{3 t_i^f}, \quad b_i = \frac{2(S_i - v_i^0 t_i^f)}{2(t_i^f)^2}. \notag
\end{align}
The value of $t_i^f$ guarantees that an unconstrained trajectory satisfies all the state, control, and safety constraints, and the boundary conditions \cite{Malikopoulos2020}.
In practice, for each CAV $j\in\mathcal{N}(t)$, the coordinator stores the optimal exit time $t_j^f$ and the corresponding coefficients $a_j, b_j, c_j, d_j$. It has been shown \cite{Malikopoulos2020} that there is no duality gap in \eqref{eq:tif}, and therefore the optimal solution can be derived in real time.

Each time a CAV $i$ enters the control zone, it communicates with the coordinator and gets access to the time trajectories of all CAVs that are inside the control zone to derive its optimal exit time $t_i^f$ from \eqref{eq:tif}.
Then, CAV $i$ transmits its four coefficients and $t_i^f$ back to the coordinator.
In the next section, we present a brief tutorial of applying our control framework to a multi-lane roundabout scenario, and we discuss several important insights that come from running scaled experiments.

\section{Tutorial: Roundabout Case Study}

\begin{figure}[htbp]
    \centering
    \includegraphics[width=0.60\textwidth]{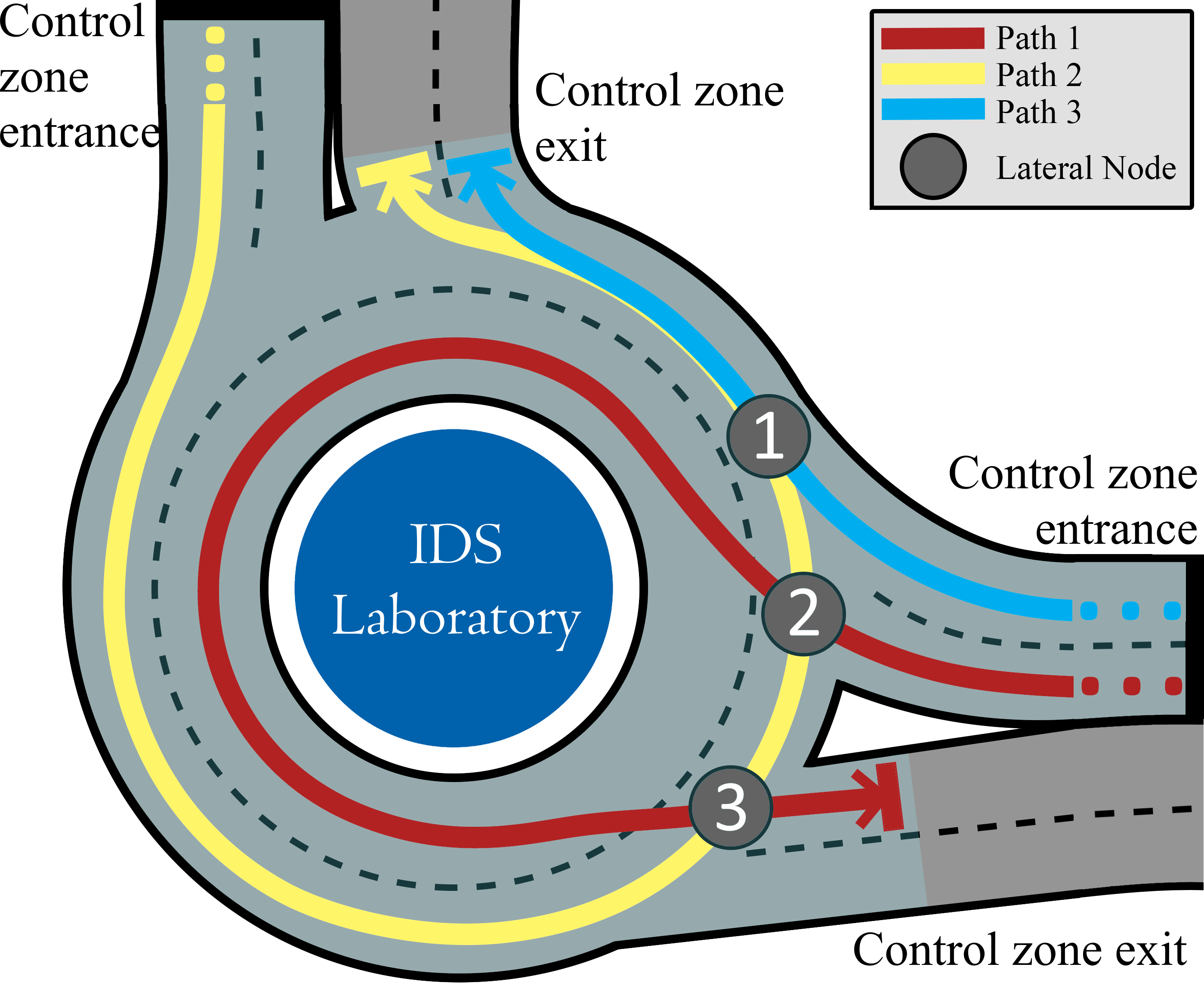}
    \caption{A schematic of the roundabout scenario. The highlighted control zone continues upstream from the roundabout at both entrances.}
    \label{fig:roundabout}
\end{figure}

To illustrate the implementation of our framework, we performed experiments in one of the two multi-lane roundabouts of IDS$^3$C (Fig. \ref{fig:roundabout}) using three CAVs per path.
Figure \ref{fig:roundabout} shows three paths with three conflict points that have a potential for lateral collisions, which we denote as lateral nodes. 
The length of the control zone for each path is $5.3$ m, $5.8$ m, and $3.8$ m ($132.5$ m, $145$ m, $95$ m scaled), respectively.
The CAVs initially operate with an intelligent driver model controller \cite{Treiber2000}, and switch to our control framework when entering the control zone.
Each CAV then determines its time trajectory by solving \eqref{eq:tif} numerically. 
The CAV follows this trajectory through the control zone. Upon exiting the control zone, it reverts to the intelligent driver model and loops back around toward the control zone entrance.
For the experiments we used the following parameters: $v_{\max} = 0.5$ m/s ($28$ mph full scale), $v_{\min} = 0.15$ m/s ($8.4$ mph full scale), $u_{\max} = 0.45$ m/s$^2$ ($11$ m/s$^2$ full scale), and $u_{\min} = - u_{\max}$.
To ensure safety, we selected a time gap of $1.0$ s and a minimum standstill distance of $0.07$ m (approx. $1$ car length). Our  framework yields an average computation time of 2.14 ms, with a maximum of 3.4 ms when a CAV plans its trajectory.
To quantify the effect of noise and disturbances acting on the system, we repeated the experiment five times.
Furthermore, we precisely timed the release of the CAVs into the roundabout such that lateral collisions would occur without intervention.
Supplementary videos of the roundabout experiment can be found at \url{https://sites.google.com/view/ud-ids-lab/csm}.

Minimum and average speed and travel time results for the five experiments are summarized in Table \ref{tab:results}.
Note that the minimum speed of all CAVs is $0.12$ m/s ($7$ mph at full scale) across all experiments using our control framework, which demonstrates that stop and go driving has been completely eliminated.
Additionally, the average speed of CAVs is $0.42$ m/s ($24$ mph at full scale), which implies that most CAVs travel near $v_{\max} = 0.5$ m/s.
The error between desired and actual exit time varies between $2-4$\%, which stems from the tracking error in the CAV's low-level controller.

\begin{table}[htbp]
    \caption{Minimum and average speed and travel time results for the $5$ experiments. The root mean square error (RMSE) of the actual exit time compared to the desired exit time from the control zone averaged over all CAVs in each experiment is provided.
    }
    \centering
    \begin{tabular}{c|ccc}
        Experiment & $v_{\min}$ [m/s] & $v_{\text{avg}}$ [m/s] & Travel Time RMSE \\
        \toprule
        1 & 0.16 & 0.41 & 2.71\,\% \\
        2 & 0.27 & 0.45 & 1.54\,\% \\
        3 & 0.18 & 0.41 & 4.03\,\% \\
        4 & 0.12 & 0.43 & 1.92\,\% \\
        5 & 0.21 & 0.42 & 1.38\,\%
    \end{tabular}
    \label{tab:results}
\end{table}

\begin{figure}[htbp]
    \centering
    \includegraphics[width=0.70\textwidth]{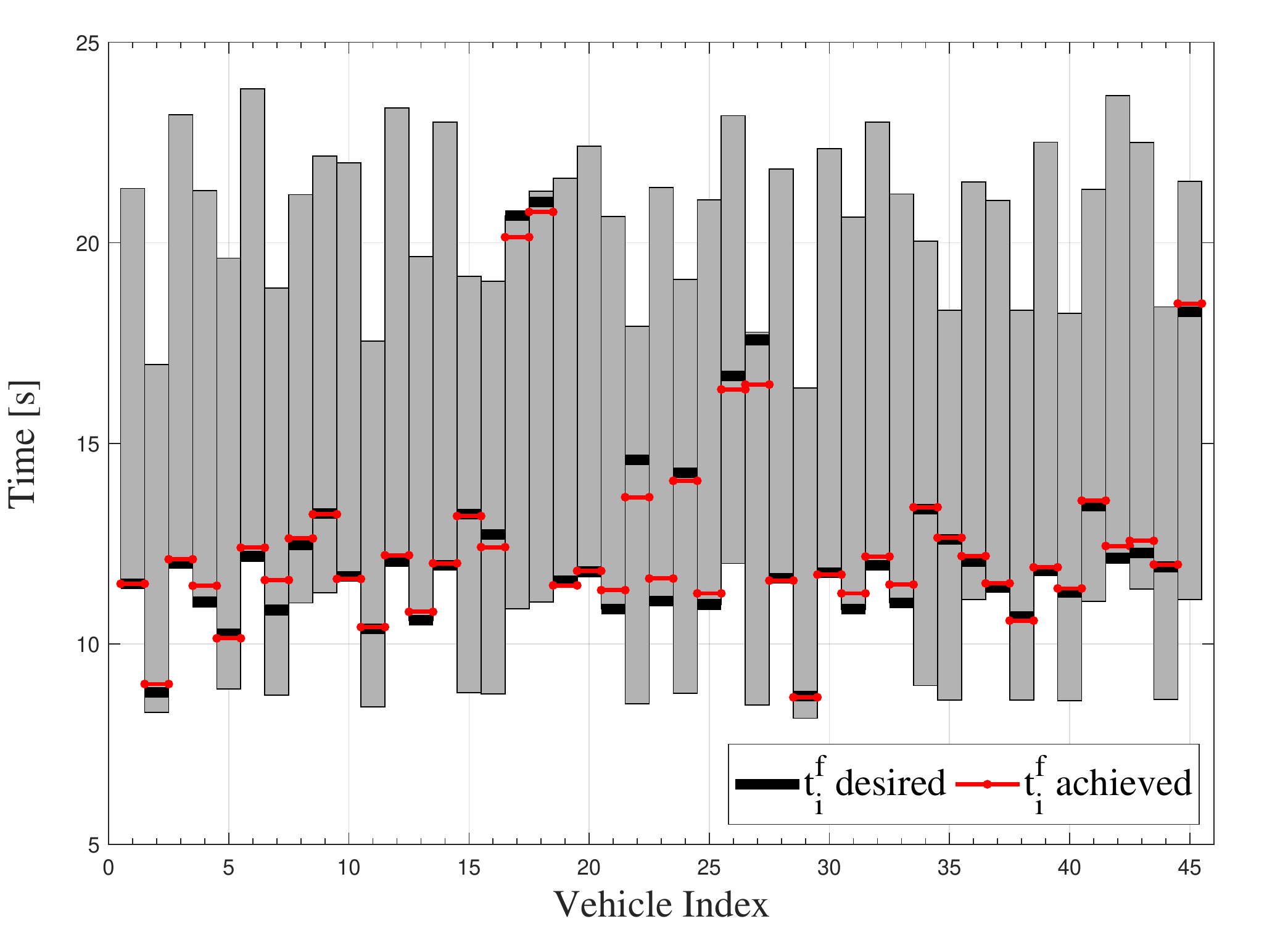}
    \caption{Planned and achieved exit time for each vehicle over all experiments. The grey bars shows the range of admissible $t_i^f$ from the state and control constraints. Every $9$ vehicles corresponds to a single experiment; they are sorted in ascending order by departure time from the control zone.}
    \label{fig:tifResult}
\end{figure}

The exit time values for each CAV are visualized in Fig. \ref{fig:tifResult}, showing the variation between the simulated and actual behavior of each CAV. The grey bars represent the feasible space of $t_i^f$, the wide black bars correspond with the planned value of $t_i^f$, and the thin red bars show the actual value of $t_i^f$ achieved by each CAV.
This effect of tracking error is visible in Table \ref{tab:results}, where the minimum achieved speed is slightly lower than the minimum speed imposed on the reference trajectory.
Figure \ref{fig:tifResult} also demonstrates how some scenarios can lead to a very small feasible space, i.e., an exit time near the maximum. This can be seen in vehicles $17, 18$, and $27$.
This motivates the introduction of a regularization zone upstream, which could influence the initial state of each CAV in the control zone to enlarge its feasible space.

\begin{figure}[htbp]
   \begin{subfigure}[htbp]{0.5\linewidth}
        \centering
        \includegraphics[width=0.95\linewidth]{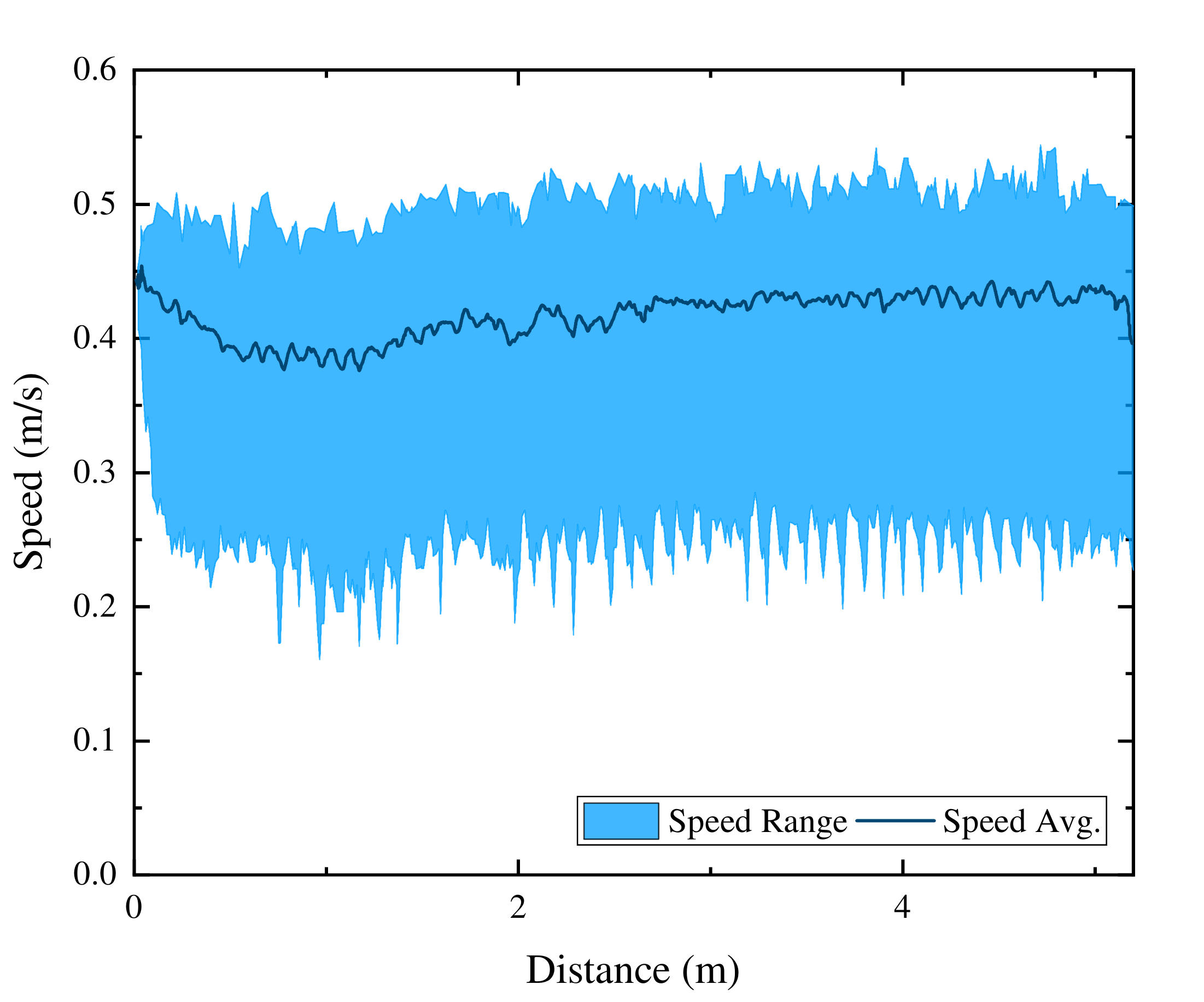}
        \caption{}
        \label{a}
   \end{subfigure}
   \begin{subfigure}[htbp]{0.5\linewidth}
        \centering
        \includegraphics[width=0.95\linewidth]{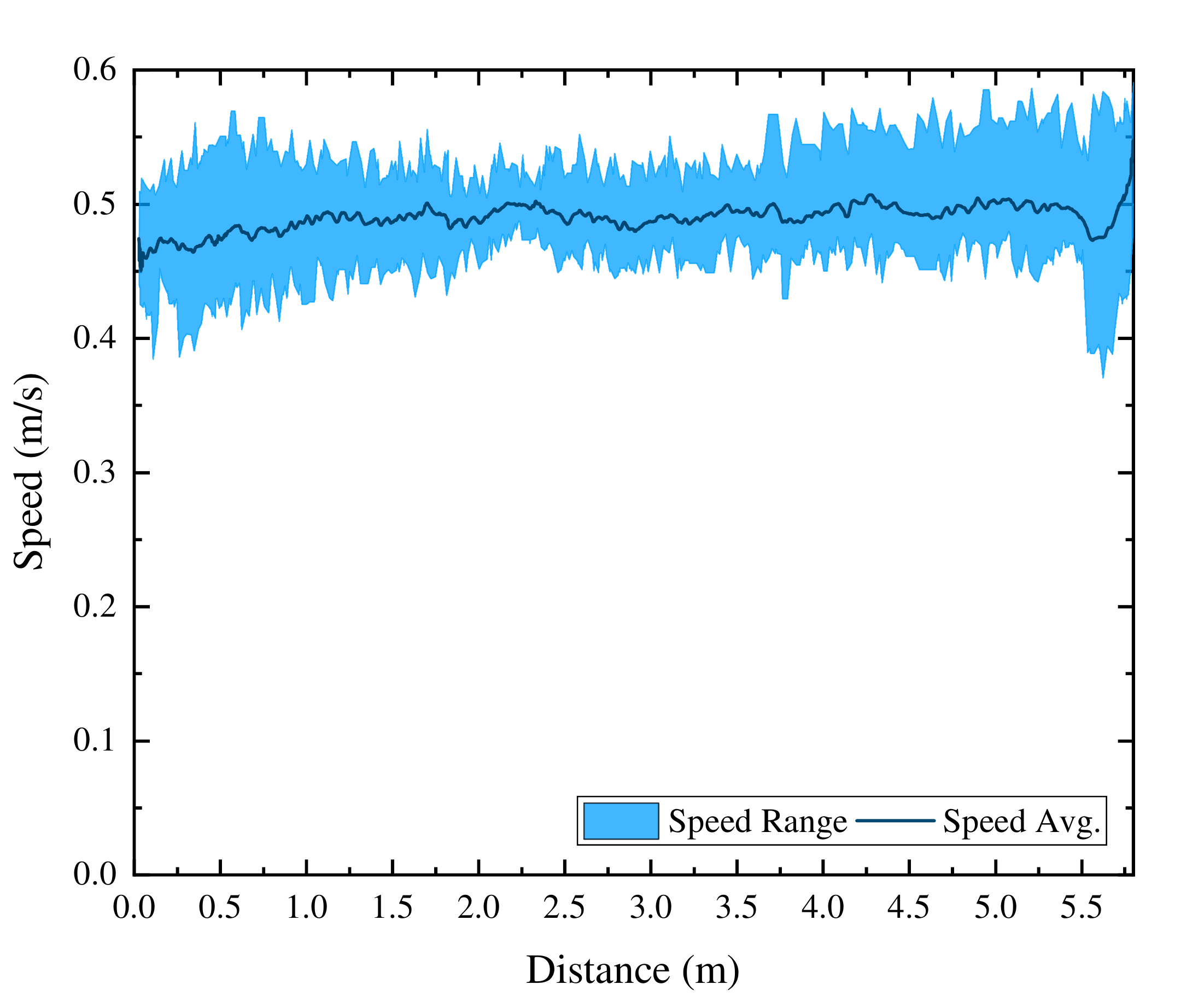}
        \caption{}
        \label{b}
   \end{subfigure}
    \begin{subfigure}[htbp]{\linewidth}
        \centering
        \includegraphics[width=0.5\linewidth]{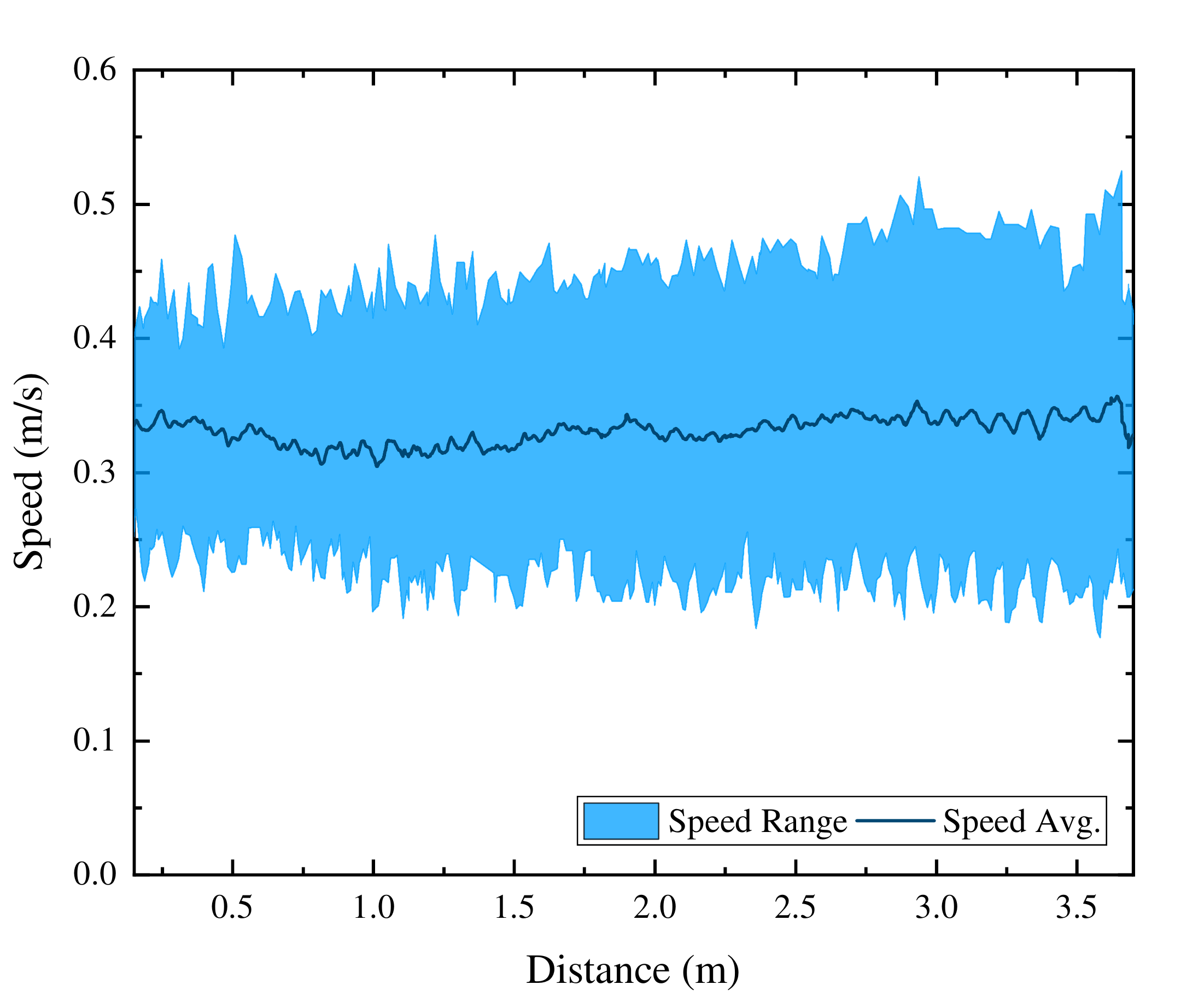}
        \subcaption{}
        \label{c}
   \end{subfigure}
    \caption{Speed range and average for all CAVs on \subref{a} path $1$, \subref{b} path $2$, and \subref{c} path $3$ across all experiments in the multi-lane roundabout.}
    \label{fig:trajectoryBlock}
\end{figure}

Finally, the average, maximum, and minimum speed for each CAV across all experiments are given in Fig. \ref{fig:trajectoryBlock}. Each figure corresponds to a single path (see Fig. \ref{fig:roundabout}) and considers $15$ CAVs ($3$ CAVs per path over five experiments). The CAVs' positions are taken directly from VICON and numerically derived using a first-order method.
From Fig. \ref{fig:trajectoryBlock}, the average speed for CAVs on each path is very close to constant.
Path $1$ shows the most variance, which is due to the distance between collision nodes $2$ and $3$ on path $1$ (see Fig. \ref{fig:roundabout}).
In order for a CAV $i\in\mathcal{N}(t)$ that is traveling along path $1$ to reduce its arrival time at node $2$, it must make a proportionally larger reduction in the value of $t_i^f$.
This is a side effect of enforcing the unconstrained trajectory on each CAV over the entire control zone.
Additionally, the entrance to the control zone along path $3$ follows a sharp right turn.
This results in relatively lower average speed in Fig. \ref{fig:trajectoryBlock}\subref{c}, as the dynamics of the CAVs reduce their speed while turning, causing them to enter the control zone at a lower initial speed.
Finally, there are instances in Fig. \ref{fig:trajectoryBlock}\subref{b} where the maximum vehicle speed surpasses the speed limit.
This is a result of stochasticity in the vehicle dynamics and sensing equipment, as well as environmental disturbances, on our deterministic controller. This analysis has motivated the development of an enhanced framework for trajectory generation of CAVs that accounts for noise, disturbances \cite{chalaki2021Reseq}, communication delay \cite{mahbub2022NHM}, and low-level tracking errors \cite{chalaki2021RobustGP,chalaki2020TCST}.

Next, we present a high-level overview and application of our control framework in a full transportation corridor in IDS$^3$C.
\section{Transportation Corridor}
\begin{figure}[htbp]
    \centering
    \includegraphics[width=0.60\textwidth]{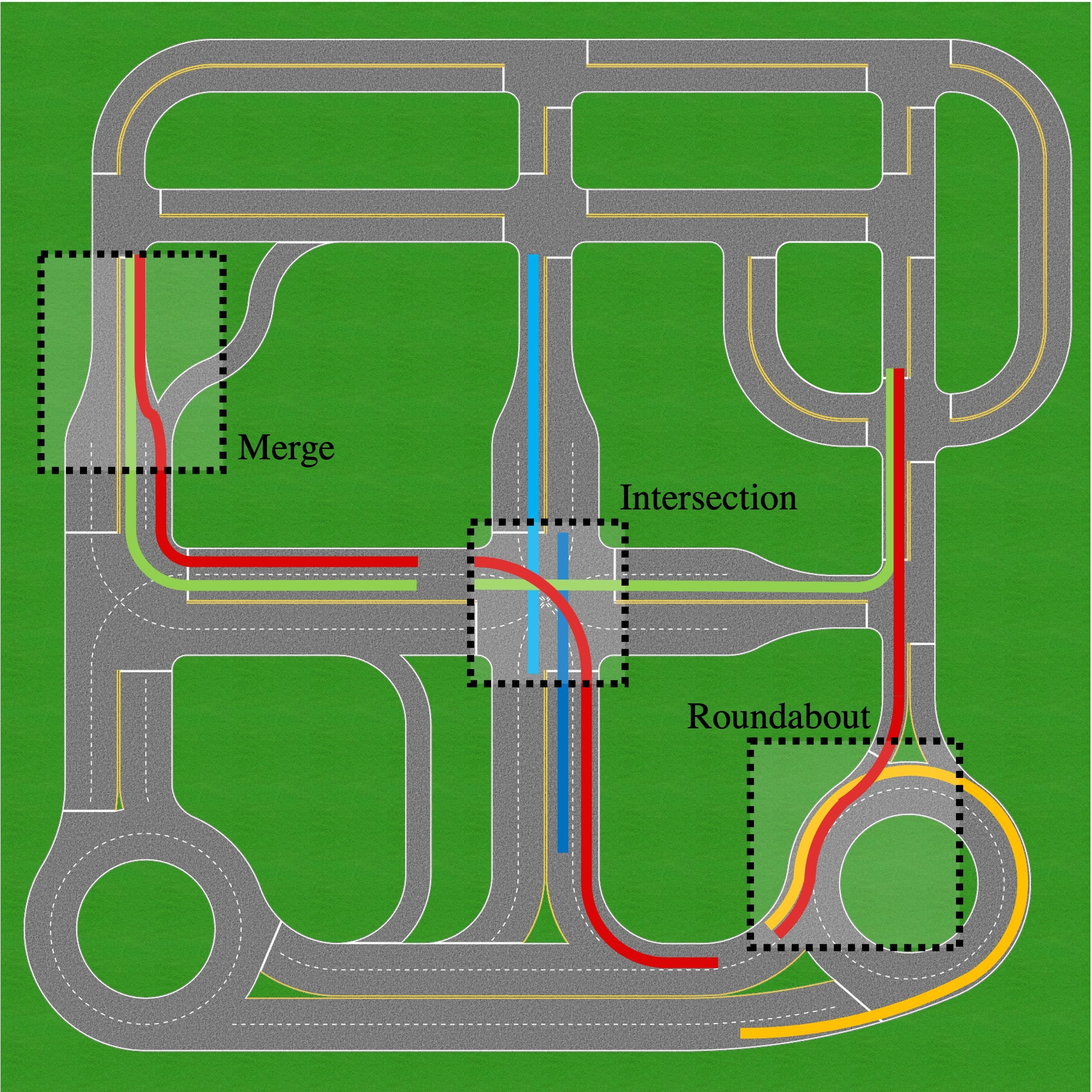}
    \caption{Corridor experiment where the ego-CAVs (red path) must navigate a roundabout, intersection, and merging roadway. The paths are only colored where they pass through a control zone, and the segments belonging to the same path have a shared color.}
    \label{fig:corridor}
\end{figure}

In this section, we apply our control framework in a transportation corridor in IDS$^3$C using a fleet of $15$ CAVs.
The corridor is shown in Fig. \ref{fig:corridor}, where $3$ ego-CAVs are released along the red path (starting in the northeast of the IDS$^3$C) and travel through a roundabout, an intersection, and a merging roadway.
At each traffic scenario, we release $3$ additional CAVs per path (as indicated in Fig. \ref{fig:corridor}) to create congestion.
The traffic scenarios were specifically selected so that upon entering the control zone, each CAV would have approximately $3$ m ($75$ m scaled) to adjust their speed before reaching a conflict point.
This also allowed us to consider each coordinator and control zone independently, as the control zone length was sufficiently long to neglect the influence of another upstream control zone.
In the baseline case, we replaced the roundabout and merging zone coordinators with yield signs.
In both scenarios, the merging vehicles yield to any vehicle within $0.4$ m of the merging zone ($10$ m scaled, approx. $4$ car lengths).
To manage the intersection in the baseline case, we implemented a four-way stop with a first-in-first-out queue. Namely, whenever a vehicle enters a line segment leading up to the intersection, it is added to the queue.
When the merging zone contains no vehicles, if the front vehicle has come to a complete stop, it is removed from the queue and allowed to pass through the merging zone.
We have taken this approach to the intersection in order to avoid any bias that may be introduced into our results by the timing of a traffic light.

Finally, to ensure a fair comparison, we set the speed limit for the entire city to $0.5$ m/s (approximately $30$ mph scaled) in both tests.
In our framework, we impose a maximum speed of $0.3$ m/s (approximately $15$ mph scaled) outside of the control zone.
This ensures that the vehicles enter the control zone at a speed lower than $v_{\max}$, and gives them the opportunity to accelerate through the control zone.

\begin{figure}[htbp]
    \centering
    \includegraphics[width=0.70\textwidth]{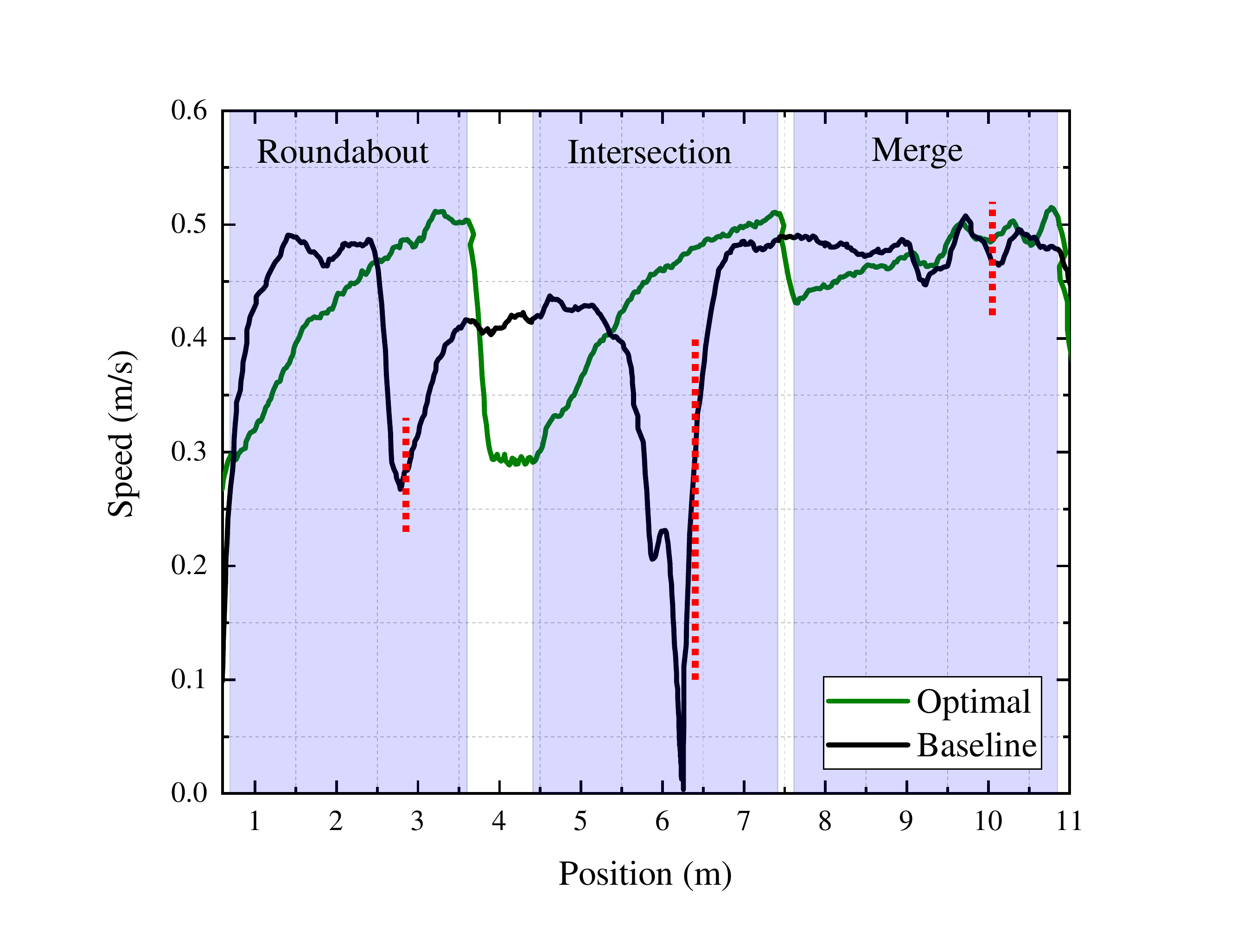}
    \caption{Speed vs position graph for the front ego vehicles in the optimal control and baseline cases. Blue highlighted areas are within each of the control zones in the optimal case, and the vertical dashed lines correspond to the location of stop and yield signs in the baseline case.
    }
    \label{fig:VvsR}
\end{figure}

\begin{figure}[htbp]
    \centering
    \includegraphics[width=0.70\textwidth]{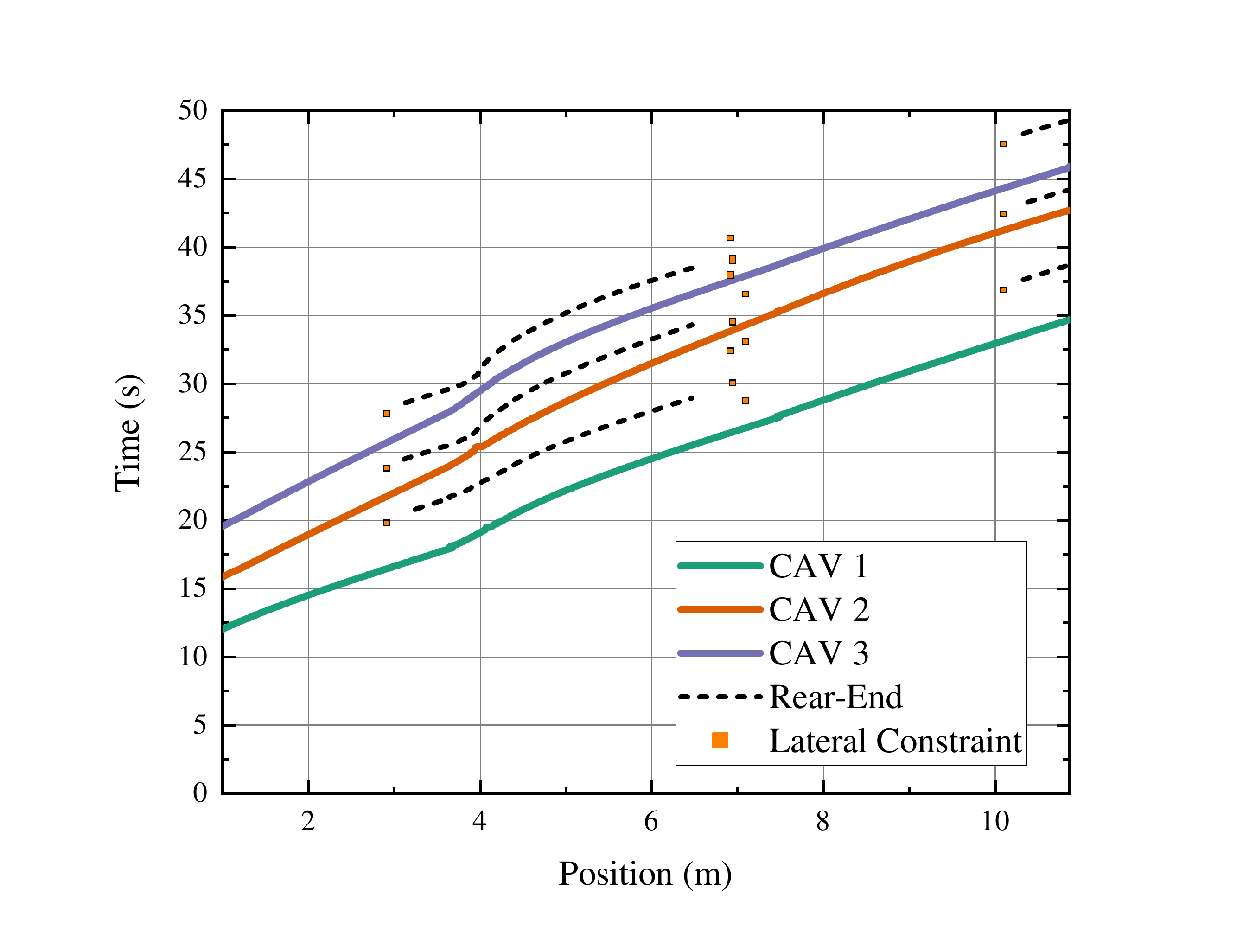}
    \caption{Time vs position graph for the ego CAVs in the optimal control case. Solid lines correspond to the CAV trajectories, dashed lines correspond to CAVs that merge onto the ego-path, and orange boxes correspond to time intervals when a lateral conflict point is occupied by another CAV.}
    \label{fig:TvsR}
\end{figure}

Figure \ref{fig:VvsR} shows that despite the apparent advantage of the baseline case's higher speed limit, the ego-CAV maintains a higher average speed in the optimal control case, and stop-and-go driving has been completely eliminated. 
Furthermore, Fig. \ref{fig:TvsR} shows that the ego-CAVs do not activate any safety constraints throughout the experiment.
Additional videos and figures of the experiment can be found at \url{https://sites.google.com/view/ud-ids-lab/csm}.

\section{Concluding Remarks}
In this article, we introduced the IDS$^3$C, a research and educational robotic scaled ($1$:$25$) testbed capable of safely validating control approaches beyond simulation in applications related to emerging mobility systems. 
This testbed can help us prove new mobility concepts and understand the implications of errors/delays in vehicle-to-vehicle and vehicle-to-infrastructure communication. IDS$^3$C can help us develop and implement control algorithms for coordinating CAVs in different traffic scenarios, such as intersections, merging roadways, speed reduction zones, roundabouts, and transportation corridors.

On the educational and outreach fronts, IDS$^3$C has been used to (a) train and educate graduate students by exposing them to a balanced mix of theory and practice, (b) integrate the research outcomes into existing courses, (c) involve undergraduate students in research, (d) create interactive educational demos, and (e) reach out to high-school students.  IDS$^3$C has been a research and educational catalyst for motivating interest in undergraduate and high-school students in science, technology, engineering, and mathematics.

We also provided an overview of a control framework for coordinating CAVs.
We demonstrated its effectiveness in IDS$^3$C in a  multi-lane roundabout using $9$ CAVs and in a corridor consisting of a roundabout, an intersection, and a merging roadway using $15$ CAVs. Ongoing research considers enhancing the framework by incorporating uncertainty originated from the vehicle surroundings  \cite{chalaki2021RobustGP,chalaki2020TCST}, and the effects of errors and delays in vehicle-to-vehicle and vehicle-to-infrastructure communication \cite{mahbub2022NHM}. Another direction of current research considers how to operate the CAVs in a way to indirectly control  human-driven vehicles and force them to form platoons led by CAVs \cite{mahbub2021_platoonMixed,Beaver2021Constraint-DrivenStudy}.

\section*{Acknowledgement(s)}

The authors would like to acknowledge Raymond Zayas and Amanda Kelly for their effort in designing, building, testing, and maintaining the newest generation of CAVs used in this paper. The authors would also like to thank Joel Diaz Goenaga for creating animations of the scaled CAV exploded view.

\bibliographystyle{Resources/IEEEtran.bst}
\bibliography{References/IDS_Publications_03102022.bib,References/Ref.bib,References/OtherTestBeds.bib,References/newRef.bib}

\newpage
\section{Author Biography}
	{Behdad Chalaki} (S’17) received the B.S. degree in mechanical engineering from the University of Tehran, Tehran, Iran in 2017. He received the M.S. degree from the department of mechanical engineering at the University of Delaware, Newark, USA, in 2021. He is currently a Ph.D. candidate in the Information and Decision Science Laboratory in the Department of mechanical engineering at the University of Delaware. His primary research interests lie at the intersections of decentralized optimal control, statistics, and machine learning, with an emphasis on transportation networks. In particular, he is motivated by problems related to improving traffic efficiency and safety in smart cities using optimization techniques. He is a student member of IEEE, SIAM, ASME, and AAAS.

{Logan Beaver} (S'17) received the B.S. degree in mechanical engineering from the Milwaukee School of Engineering, Milwaukee, WI, in 2015 and the M.S. degree in mechanical engineering from Marquette University,
Milwaukee, WI, in 2017.
He is currently a Ph.D. candidate in the Information and Decision Science Laboratory in the department of mechanical engineering at the University of Delaware.
His research focuses on the detection, stabilization, and control of emergent behavior in multi-agent and swarm systems.
He is a graduate student member of IEEE, SIAM, ASME, and AAAS.

{A M Ishtiaque Mahbub} (S’17) received the B.S. degree in mechanical engineering from Bangladesh University of Engineering and Technology, Bangladesh in 2013 and M.Sc. in computational mechanics from University of Stuttgart, Germany in 2016. He is currently pursuing a Ph.D. degree in mechanical engineering at the Information Decision Science Laboratory at University of Delaware under the supervision of Prof. Andreas A. Malikopoulos. His research interests include, but are not limited to, optimization and control with an emphasis on applications related to connected automated vehicles, hybrid electric vehicles, and intelligent transportation systems. He has conducted several internships at National Renewable Energy Lab (NREL), Robert Bosch LLC (USA), Robert Bosch GmbH (Germany) and Fraunhofer IPA (Germany). He is a student member of IEEE, SIAM and SAE.

{Heeseung Bang} (S'18) received the B.S. degree in mechanical engineering and M.S. degree in electrical engineering from Inha University, Incheon, South Korea, in 2018 and 2020. He is currently a Ph.D. student in the department of mechanical engineering at the University of Delaware. His research interests are in optimization and control of multi-agent systems, with an emphasis on autonomous mobility-on-demand systems and future mobility systems. 

{Andreas A. Malikopoulos}
	(S'06--M'09--SM'17) received the Diploma in mechanical engineering from the National Technical University of Athens, Greece, in 2000. He received M.S. and Ph.D. degrees from the department of mechanical engineering at the University of Michigan, Ann Arbor, Michigan, USA, in 2004 and 2008, respectively. 
	He is the Terri Connor Kelly and John Kelly Career Development Associate Professor in the Department of Mechanical Engineering at the University of Delaware (UD), the Director of the Information and Decision Science (IDS) Laboratory, and the Director of the Sociotechnical Systems Center. Before he joined UD, he was the Deputy Director and the Lead of the Sustainable Mobility Theme of the Urban Dynamics Institute at Oak Ridge National Laboratory, and a Senior Researcher with General Motors Global Research \& Development. His research spans several fields, including analysis, optimization, and control of cyber-physical systems; decentralized systems; stochastic scheduling and resource allocation problems; and learning in complex systems. The emphasis is on applications related to smart cities, emerging mobility systems, and sociotechnical systems. He has been an Associate Editor of the IEEE Transactions on Intelligent Vehicles and IEEE Transactions on Intelligent Transportation Systems from 2017 through 2020. He is currently an Associate Editor of Automatica and IEEE Transactions on Automatic Control. He is a member of SIAM, AAAS, and a Fellow of the ASME.

\end{document}